\documentclass[a4paper,10pt]{article}
\usepackage{amsmath,amsthm,amssymb}
\usepackage{graphicx,subfigure}
\usepackage{url}

\usepackage{algorithm2e}
\usepackage{listings}
\usepackage{color}

\usepackage{multicol}
\def\twoplot[#1]#2#3#4#5{
\begin{figure}[h]
\begin{multicols}{2}
\begin{center}
    \includegraphics*[#1]{#2}
    \caption{\label{#2} #4}
\end{center}
\begin{center}
    \includegraphics*[#1]{#3}
    \caption{\label{#3} #5}
\end{center}
\end{multicols}
\end{figure}
}

\usepackage[colorlinks,bookmarksopen,bookmarksnumbered,citecolor=red,urlcolor=red]{hyperref}
\hypersetup{pdftitle={FullSWOF\_Paral: Comparison of two parallelization strategies (MPI and SKELGIS) on a software designed for
 hydrology applications},
bookmarks=true,
pdftoolbar=true,
pdfmenubar=true,
pdfauthor={S. Cordier, H. Coullon, O. Delestre, C. Laguerre, M.-H. Le, D. Pierre, G. Sadaka},
pdfsubject={applied mathematics, hydraulics, parallelization},
pdfcreator={Delestre},
pdfproducer={Delestre},
pdfkeywords={Shallow-Water equation} {Saint-Venant} {parallelization} {MPI} {SKELGIS} {skeleton} {domain decomposition}
 {Malpasset} {finite volume} {well-balanced} {hydrostatic reconstruction} {HLL numerical flux} {MUSCL reconstruction}}

\title{FullSWOF\_Paral: Comparison of two parallelization strategies (MPI and SKELGIS)\\ on a software designed for hydrology
 applications\footnote{Authors thank CaSciModOT, FED 4222 University of Orl\'eans and AMIES  \url{http://www.agence-maths-entreprises.fr/}.}}

\author{{S. Cordier}\footnote{MAPMO, f\'ed\'eration Denis Poisson, University of Orl\'eans, France, e-mail :
 stephane.cordier@math.cnrs.fr}, {H. Coullon}\footnote{LIFO, University of Orl\'eans \& G\'eo-Hyd, France, e-mail :
 helene.coullon@univ-orleans.fr}, O. Delestre\footnote{Lab. J.A. Dieudonn\'e \& EPU Nice Sophia, University of
 Nice, France, e-mail : delestre@math.unice.fr}, {C. Laguerre}\footnote{MAPMO, f\'ed\'eration Denis Poisson, University of Orl\'eans,
 France, e-mail : christian.laguerre@math.cnrs.fr},\\ 
 {M.-H. Le}\footnote{MAPMO \& BRGM, University of Orl\'eans, France, e-mail : Minh.Hoang.Le@math.cnrs.fr},
 D. Pierre\footnote{G\'eo-Hyd, 101 rue Jacques Charles, 45160 Olivet, France, e-mail : daniel.pierre@geo-hyd.com}
 \; and G. Sadaka\footnote{LAMFA, University of Picardie, France, e-mail : georges.sadaka@u-picardie.fr}}

\begin{document}
\maketitle

\begin{abstract} In this paper, we perform a comparison of two approaches for the parallelization of an existing, free software, FullSWOF\_2D
 (\url{http ://www.univ-orleans.fr/mapmo/soft/FullSWOF/} that solves shallow water equations for applications in hydrology) based on a domain
 decomposition strategy. The first approach is based on the classical MPI library while the second approach uses Parallel Algorithmic Skeletons and
 more precisely a library named SkelGIS (Skeletons for Geographical Information Systems). The first results presented in this article show that
 the two approaches are similar in terms of performance and scalability. The two implementation strategies are however very different
 and we discuss the advantages of each one.\end{abstract}

\section{Problematics}
\label{sect:problematics}
We are interested in overland flow simulations. For this kind of flow simulation, several methods are used from empirical models to
 physically based models. Two physical models are often used to model overland flow kinematic (KW) and diffusive wave (DW) equations
 \cite{Moussa00}, \cite{Novak10}. But following \cite{Zhang89}, \cite{Esteves00}, \cite{Fiedler00}, we choose to use the shallow water
 (Saint-Venant \cite{saintvenant71}) physical model. Indeed KW and DW models may give poor results in terms of water heights and velocities
 in case of mixed subcritical and supercritical flows. Despite the computational cost, shallow water model is mandatory. MacCormack
 scheme is widely used to solve shallow water (SW) equations \cite{Zhang89}, \cite{Esteves00}, \cite{Fiedler00}. But it neither guarantees the
 positivity of water depths at the wet/dry transitions, nor preserves steady states (not well-balanced \cite{Greenberg96}) as noticed in 
 \cite{Rousseau}. In industrial codes (ISIS, Canoe, HEC-RAS, MIKE11, ...), SW equations are often solved under non-conservative
 form \cite{Cunge64,Cunge80,Novak10} with either Preissmann scheme or Abbott-Ionescu scheme. Thus transcritical flows and hydraulic jumps are not
 solved properly. In order to cope with all these problems, we have developed FullSWOF\_2D (Full Shallow Water equations for Overland Flow) an
 object-oriented C++ code (free software and GPL-compatible license CeCILL-V2\footnote{\url{http://www.cecill.info/index.en.html}}) in the
 framework of the multidisciplinary project METHODE (see \cite{Delestre10b} and \url{http://www.univ-orleans.fr/mapmo/methode/}). The source code
 is available at \url{http://www.univ-orleans.fr/mapmo/soft/FullSWOF/}. This software is based on the physical model of shallow water equations
 (Saint-Venant system \cite{saintvenant71}) for the water runoff, coupled with Green-Ampt's infiltration model \cite{Green11}. The set of shallow
 water equations is solved thanks to a well-balanced finite volume method. This enables to catch and to preserve steady states such as puddles
 and lake at rest. Moreover the method preserves the positivity of the water height. Validations of FullSWOF\_2D have already been performed
 on analytical benchmarks (SWASHES \cite{Delestre2012}), on experimental data and on real events at small scales
 (parcels \cite{DBLP:journals/corr/abs-1204-3210,Delestre10b}). We now aim at simulating at bigger scales such as watershed, river valley or at
 small scales with fine details. Thus the computing time becomes bigger and the code needs to be parallelized. We want to compare
 two approaches: the first one is a "classical" approach based on a master-slave architecture using MPI and the second one uses skeletal algorithms
 (SkelGIS, Skeletons for Geographical Information Systems, which is under development by H\'el\`ene Coullon, a CIFRE PhD student with G\'eo-Hyd
 company). These approaches will be compared both in terms of performance and scalability and the advantages
 of each choice will be briefly discussed.

\section{Physical model}
\label{sect:model}
\subsection{General settings}

\begin{figure}[htbp]
\begin{center}
\includegraphics[width=0.9\textwidth]{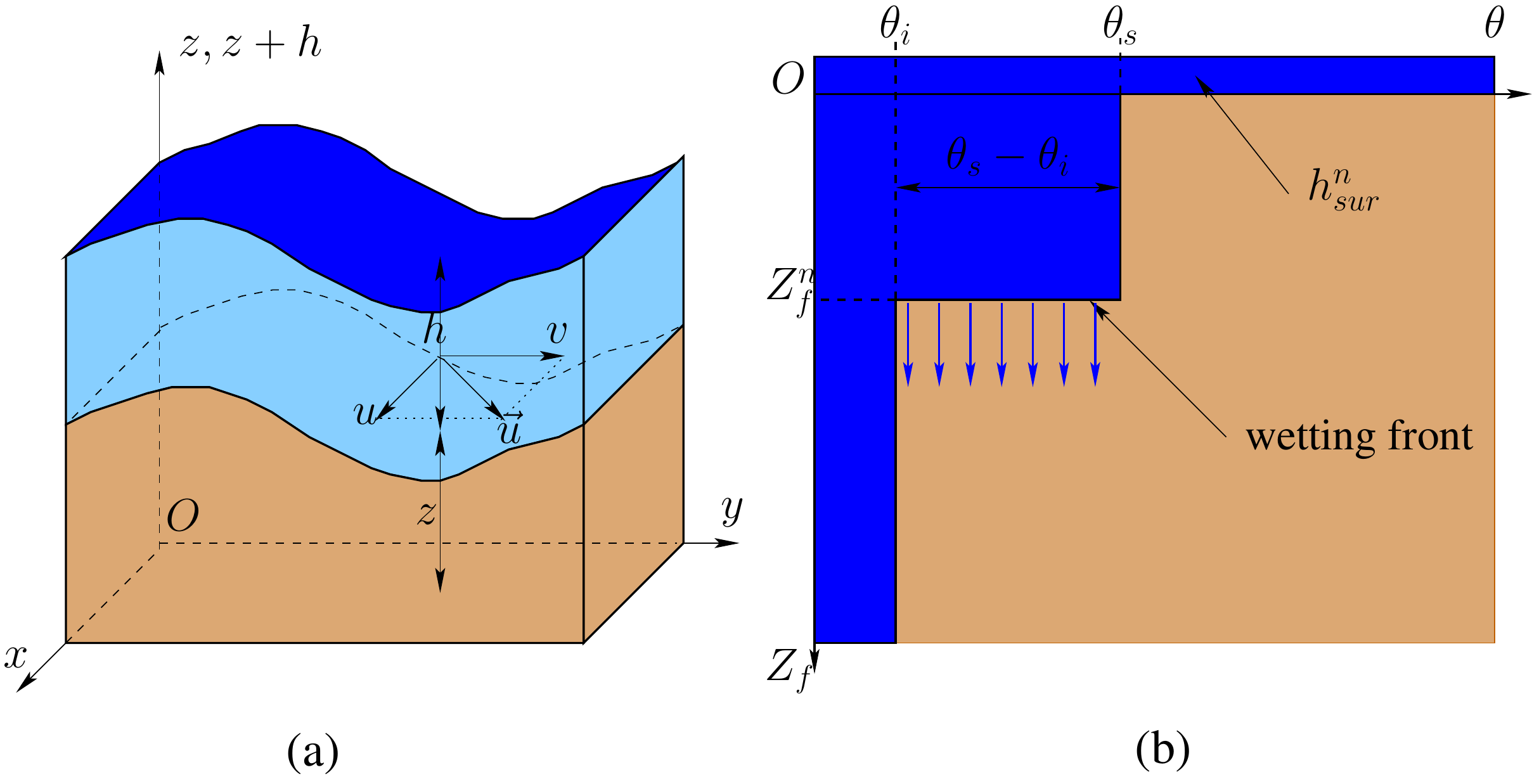}
\caption{Variables of (a) 2D shallow water equations (SW2D) and (b) Green-Ampt infiltration model.}
\label{SW2D-fig}
\end{center}
\end{figure}

As in \cite{Esteves00,Fiedler00}, we model overland flow thanks to the 2D shallow water equations (SW2D). Shallow water equations
 (or Saint-Venant system) have been proposed by Adh\'emar Barr\'e de Saint-Venant in 1871 
in order to model flows in channels \cite{saintvenant71} in one dimension in space. 
Nowadays, they are used to model flows in various contexts, such as: overland flow \cite{Esteves00,Tatard08},
rivers \cite{Goutal02,Burguete04}, flooding \cite{Caleffi03,Delestre09}, dam breaks \cite{Alcrudo99,Valiani02}, 
nearshore \cite{Borthwick01,Marche05}, tsunami \cite{George06,Kim07,Popinet11}.
These equations consist in a nonlinear system of partial differential equations (PDE-s), more precisely
conservation laws describing the evolution of the height ($h(t,x,y)$ [L]) and the horizontal components of the vertically averaged velocity of the
 fluid ($\vec{u}(t,x,y)=(u(t,x,y),v(t,x,y))^t$ [L/T]) as illustrated on figure \ref{SW2D-fig}-a. This complete set of conservation laws writes
\begin{equation}
\left\{
\begin{array}{l}
 \partial_t h+\partial_x (hu) + \partial_y (hv)= R-I\\
\partial_t(hu)+\partial_x\left( hu^2+gh^2/2\right)+\partial_y \left(huv \right)=gh\left({S_0}_x-{S_f}_x \right)\\
\partial_t(hv)+\partial_x\left(huv\right)+\partial_y \left( hv^2+gh^2/2 \right)=gh\left({S_0}_y-{S_f}_y \right)\label{SW2D}
\end{array}
\right.
\end{equation}
where 
\begin{itemize}
 \item $g=9.81\,\text{m}/\text{s}^2$ is the gravity constant; 
 \item $z(x,y)$ is the topography [L] and we denote by ${S_0}_x=-\partial_x z(x,y)$ (resp. ${S_0}_y=-\partial_y z(x,y)$) the opposite of the
 slope in the $x$ (resp. $y$) direction. Erosion is not considered here, so the topography is a fixed function of space. Equations might be added
 to SW2D model to take erosion effect into account. We get systems such as Saint-Venant Exner and Hairsine \& Rose models (for more details see
 \cite{Le12});
 \item $R(t,x,y)\geq 0$ [L/T] is the rain intensity;
 \item $I(t,x,y)$ [L/T] is the infiltration rate. It is given by another model (such as Green-Ampt \cite{Green11},
 Richards \cite{Richards31}, ...);
 \item and $\vec{S}_f=\left({S_f}_x,{S_f}_y\right)$ is the friction force vector. It may take several forms, depending on soil and flow properties.
 In hydrological and hydraulics models, two families of friction laws are mainly encountered. They are based on empirical considerations.
 On one hand, we have the family of Manning-Strickler's friction laws
 \begin{equation*}
  \vec{S}_f=C_f\dfrac{\sqrt{u^2+v^2}}{h^{4/3}}\vec{u},
 \end{equation*}
$C_f=n^2$ (resp. $C_f=1/K^2$), where $n$ (resp. $K$) is the Manning's coefficient [L$^{-1/3}$T] (resp. Strickler's coefficient [L$^{-1/3}$T]).
 On the other hand, we have the law of Darcy-Weisbach's and Ch\'ezy's family
\begin{equation*}
 \vec{S}_f=C_f \dfrac{\sqrt{u^2+v^2}}{h}\vec{u},
\end{equation*}
$C_f=f/(8g)$ (resp. $C_f=1/C^2$), where $f$ (resp. $C$) is the dimensionless Darcy-Weisbach's coefficient (resp. Ch\'ezy's coefficient
 [L$^{1/2}$/T]). The friction may depend on the space variable, especially on large parcels and watersheds. Values of friction coefficients
 depend on the considered type of soil and are tabulated in references such as \cite{Chow59}.
\end{itemize}

System \eqref{SW2D} will be solved on a Cartesian grid thanks to a method of lines. Thus the numerical strategy consists in choosing a numerical
 method adapted to the properties of the shallow water system in 1D (SW1D). Then the generalization to SW2D is straightforward. So in next sections,
 we will consider SW1D system.

\subsection{Hyperbolicity}

As previously argumented, we place ourselves in the one-dimensional case. SW1D system writes
\begin{equation}
\left\{
\begin{array}{l}
 \partial_t h+\partial_x (hu)= R-I\\
\partial_t(hu)+\partial_x\left( hu^2+gh^2/2\right)=gh\left({S_0}-{S_f} \right)\label{SW1D}
\end{array}
\right..
\end{equation}

The left-hand side of system \eqref{SW1D} is a transport operator, corresponding to the flow of an ideal fluid in a flat channel, without friction,
 rain or infiltration. This corresponds to the model introduced by Saint-Venant \cite{saintvenant71} which contains several flow properties.
 To emphasize these properties, we rewrite the one-dimensional homogeneous system under vectors form
\begin{equation}
 \partial_t U+\partial_x F(U)=0,\;\text{where}\;U=
\left(\begin{array}{c}
       h\\
      hu
      \end{array}\right)=
\left(\begin{array}{c}
       h\\
      q
      \end{array}\right),\;
F(U)=\left(\begin{array}{c}
            hu\\
        hu^2+gh^2/2
           \end{array}\right),\label{SW1DH}
\end{equation}
with $F(U)$ the flux of the equation and $q(t,x)$ [L$^2$/T] the discharge by unit of width. The transport property of system \eqref{SW1DH} is
 clearer under the following nonconservative forms
\begin{equation*}
 \partial_t U+A(U)\partial_x U=0,\; A(U)=F'(U)
=\left(\begin{array}{cc}
        0 & 1 \\
    gh-u^2 & 2u
       \end{array}\right),
\end{equation*}
where the Jacobian matrix $A(U)=F'(U)$ is the matrix of transport coefficients. When $h>0$, matrix $A(U)$ is diagonalizable and its eigenvalues
 are
\begin{equation}
 \lambda_1(U)=u-\sqrt{gh}<u+\sqrt{gh}=\lambda_2(U).\label{eingval}
\end{equation}
System \eqref{SW1DH} is strictly hyperbolic. The eigenvalues are the velocities of the surface waves and thus are basic characteristics of
 flows. Notice that these eigenvalues collapse in one when $h=0\,\text{m}$ (for dry zones). In that case, system \eqref{SW1DH} is no longer
 hyperbolic which is difficult to deal with both from theoretical and numerical levels. Getting a numerical scheme that preserves the water
 height positivity is a necessity.

Eigenvalues \eqref{eingval} of system \eqref{SW1DH} allow to make a classification of flows, based on the relative values of the velocities of the
 waves $c=\sqrt{gh}$ and of the fluid $u$. Indeed if $|u|<\sqrt{gh}$, the characteristic velocities have opposite signs and flow informations
 propagate both upward and downward. The flow is said to be subcritical or fluvial. When $|u|>\sqrt{gh}$, the flow is supercritical or torrential
 and all informations propagate downstream. This has to be taken into account both for the numerical methods (upwinding) and the boundary
 conditions.


Numerically, boundary conditions are expected to provide numerical fluxes at the boundaries of domain. These are required in order to update $(h,q)$ on the extreme cells to the next time level. The boundary conditions may result in direct prescription of the numerical fluxes at the boundaries. Alternatively, we may prescribe the values of $(h,q)$ on the ghost cells. In this way, the  Riemann problems at the boundaries are solved and the corresponding fluxes are computed as done for the interior cells. The values of $(h,q)$ on ghost cells can be computed in function of flow regime thanks to the Riemann invariants (which remain constant along the corresponding characteristic line). For subcritical flow, the characteristic $dx/dt = u-c$ leaves the domain while the characteristic $dx/dt = u+c$ enters the domain. Thus for numerical simulations, we impose one variable ($h$ or $q$) for fluvial inflow/ouflow. In contrary, both variable must be prescribed in case of supercritical inflow whereas free boundary conditions are considered for supercritical outflow. We refer to \cite{Bristeau01} for more details.

With the presence of the source terms, an other main property has to be considered: the occurrence of steady states or stationary solutions.
 These particular flows are studied in next section.

\subsection{Steady flows}

Steady states solutions correspond to stationary flows, {\it i.e.} solutions that satisfy
\begin{equation*}
 \partial_t h=\partial_t u=\partial_t q=0,
\end{equation*}
system \eqref{SW1D} reduces in
\begin{equation}
\left\{
\begin{array}{l}
 \partial_x (hu)= R-I\\
\partial_x\left( hu^2+gh^2/2\right)=gh\left({S_0}-{S_f} \right)\label{SW1Dstat}
\end{array}
\right..
\end{equation}
To our knowledge there is no scheme designed to preserve these stationary flows. If we consider no rain $R=0$, no infiltration $I=0$ and
 no friction $S_f=0$, system \eqref{SW1Dstat} reduces in
\begin{equation*}
\left\{
\begin{array}{l}
 \partial_x (hu)= 0\\
\partial_x\left( hu^2+gh^2/2\right)=-gh\partial_x z
\end{array}
\right..
\end{equation*}
We recover Bernoulli's law
\begin{equation}
\left\{
\begin{array}{l}
 hu= q=Cst\\
\dfrac{q^2}{2gh^2}+h+z=Cst
\end{array}
\right.,\label{eq-bernoulli}
\end{equation}
where $h+z$ is the free surface water level.

Several schemes have been designed to preserve steady states \eqref{eq-bernoulli}. In general, these methods are costly in term of
 implementation and computation because, for example, they lead to solve a third order polynomial and the appropriate root must be selected
 depending on the type of flow (see {\it e.g.} \cite{Castro07}). Moreover, the preserving of steady states \eqref{eq-bernoulli} becomes more
 complicated when we require also the positivity preserving of water depths at numerical level. As the last property is a compulsory property
 for the applications we are interested in (overland flow simulations), we limit ourselves to numerical schemes that preserve the particular
 steady state corresponding to equilibrium of lake at rest
\begin{equation*}
 u=q=0 \;\text{and}\; h+z=Cst.
\end{equation*}
In that case, we have hydrostatic balance between the hydrostatic pressure and the gravitational acceleration down
 the inclined bottom $z$.

 In 1994, Berm\'udez and V\'azquez \cite{Bermudez94} are the first to identify the difficulty to preserve this steady state.
 Schemes preserving this equilibrium have the Conservation property or C-property (introduced in \cite{Bermudez94}). They obtained such a
 numerical method by modifying the Roe scheme. The topography source term is upwinded thanks to a projection on the eigenvalues of Jacobian
 matrix of the flux. Since \cite{Greenberg96}, schemes which preserve exactly at least the hydrostatic equilibrium at a discrete level are
 called well-balanced schemes. We can find a lot of well-balanced schemes for SW1D in the literature \cite{Bermudez94,Bermudez98,Greenberg96,
LeVeque98,Jin01a,Kurganov02,Gallouet03,Audusse04c,Audusse05,Noelle07,Liang09,Berthon12}. 

In case of rainfall overland flows, we have the occurrence
 of wet/dry transitions and small water heights.
 To simulate that kind of events, we need a robust and positive well-balanced scheme. Thus we have chosen a finite volume scheme based on the
 hydrostatic reconstruction (introduced in \cite{Audusse04c} and \cite{Bouchut04}), that we will detail in next section.

\section{Numerical method}
\label{sect:method}
As argued in the previous section, the numerical method will be presented in one dimension (extension to 2D being straightforward on
 Cartesian grids thanks to the method of lines). 

\subsection{Convective step}\label{sec:conv}

A finite volume discretization of SW1D, \eqref{SW1D}, writes
\begin{equation}
 U_i^*=U_i^n-\dfrac{\Delta t}{\Delta x}\left[ F_{i+1/2L}^n-F_{i-1/2R}^n-Fc_i^n\right]\label{SW1Dscheme}
\end{equation}
with $\Delta x$ (resp. $\Delta t$) the space (resp. time) step and
\begin{equation*}
 \begin{array}{l}
  F_{i+1/2L}^n=F_{i+1/2}^n+S_{i+1/2L}^n\\
  F_{i-1/2R}^n=F_{i-1/2}^n+S_{i-1/2R}^n
 \end{array},
\end{equation*}
the left and right modifications of the numerical flux $\cal F$ for the homogeneous problem (see section \ref{sec:flux})
\begin{equation*}
 F_{i+1/2}^n={\cal F}\left(U_{i+1/2L}^n,U_{i+1/2R}^n\right).
\end{equation*}
The values $U_{i+1/2L}$ and $U_{i+1/2R}$ are obtained thanks to two consecutive reconstructions. Firstly a MUSCL reconstruction
 \cite{vanLeer79,Bouchut04,Delestre10b} is performed on $u$, $h$ and $h+z$ in order to get a second order scheme in space (see section
 \ref{sec:muscl}). This gives us the reconstructed values $\left(U_-,z_-\right)$ and $\left(U_+,z_+\right)$. Secondly, we apply the hydrostatic
 reconstruction \cite{Audusse04c,Bouchut04} on the water height which allows us to get a positive preserving well-balanced scheme
\begin{equation*}
 \left\{\begin{array}{l}
         h_{i+1/2L}=\max \left(h_{i+1/2-}+ z_{i+1/2-}-\max\left(z_{i+1/2-},z_{i+1/2+}\right),0\right)\\
    U_{i+1/2L}=\left(h_{i+1/2L},h_{i+1/2L}u_{i+1/2-}\right)^t\\
         h_{i+1/2R}=\max \left(h_{i+1/2+}+ z_{i+1/2+}-\max\left(z_{i+1/2-},z_{i+1/2+}\right),0\right)\\
    U_{i+1/2R}=\left(h_{i+1/2R},h_{i+1/2R}u_{i+1/2+}\right)^t
        \end{array}\right..
\end{equation*}
We introduce
\begin{equation*}
 S_{i+1/2L}^n=\left(\begin{array}{c}
                     0\\
    	    g \left(h_{i+1/2-}^2 - h_{i+1/2L}^2\right)/2
                    \end{array}\right), \;
 S_{i-1/2R}^n=\left(\begin{array}{c}
                     0\\
		    g \left(h_{i-1/2+}^2 - h_{i-1/2R}^2\right)/2
                    \end{array}\right)
\end{equation*}
and a centered source term is added to preserve consistency and well-balancing \cite{Audusse04c,Bouchut04}
\begin{equation*}
 Fc_i^n=\left(\begin{array}{c}
             0\\
	  -g \dfrac{h_{i-1/2+}+h_{i+1/2-}}{2}\left(z_{i+1/2-}-z_{i-1/2+}\right)
            \end{array}\right).
\end{equation*}
We have to insist on the positivity and the robustness of this method. The rain source term is treated explicitly and the infiltration rate is
 obtained thanks to the Green-Ampt model \cite{Green11} (see section \ref{GreenAmpt}).

\subsection{Friction treatment}\label{sec:fric}

In this step, the friction term is taken into account with the following system
\begin{equation*}
 \partial_t U=\left( \begin{array}{c}
                      0\\
		    -gh S_f
                     \end{array}\right).
\end{equation*}
This system is solved thanks to a semi-implicit method (as in \cite{Bristeau01,Fiedler00}), which writes for the Darcy-Weisbach's law
\begin{equation*}
 h^{n+1}=h^*\;\text{and}\; q^{n+1}=\dfrac{q^*}{1+\Delta t \dfrac{f}{8}\dfrac{|q^n|}{h^n h^{n+1}}},
\end{equation*}
where $h^*$, $q^*$ and $u^*$ are the variables from the convective step. This method allows to preserve stability (under a classical CFL
 condition) and steady states at rest. Finally, these two steps are combined in a second order TVD Runge Kutta method which is
 the Heun's predictor-corrector method \cite{Shu88}. It writes
\begin{equation*}
 \begin{array}{c}
  U^*=U^n+\Delta t \Phi\left(U^n \right)\\
 U^{**}=U^*+ \Delta t \Phi \left(U^* \right)\\
U^{n+1}=\dfrac{U^n + U^{**}}{2}
 \end{array}
\end{equation*}
where $\Phi$ is the right part of \eqref{SW1Dscheme}.

\subsection{Numerical flux}\label{sec:flux}
About the homogeneous flux ${\cal F}\left(U_{i+1/2L}^n,U_{i+1/2R}^n\right)$, we can use any consistent numerical flux, for example the one
 of Godunov, Rusanov, HLL, Roe or the one obtained by the kinetic or the relaxation methods. In this work, we adopted the Harten Lax van
 Leer (HLL) flux \cite{Harten83,Bouchut04,Delestre10b} which is known to be a simple and efficient solver for both accuracy and implementation
 aspects (see~\cite{Delestre10b})
\begin{equation*}
 {\cal F}\left(U_L,U_R \right)=
\left\{\begin{array}{ll}
        F(U_L) & \text{if}\; 0\leq c_1\\
\dfrac{c_2 F(U_L)-c_1 F(U_R)}{c_2 -c_1}+\dfrac{c_1 c_2}{c_2 -c_1}\left(U_R - U_L\right) & \text{if}\; c_1<0<c_2 \\
        F(U_R) & \text{if}\; c_2 \leq 0
       \end{array}\right.,
\end{equation*}
with two parameters $c_1<c_2$ 
which are the approximations of slowest and fastest wave speeds respectively. We refer to \cite{Batten1997} for further discussion on the wave
 speed estimates. In this paper, we use
\begin{equation*}
c_{1}={\inf\limits_{U=U_L,U_R}}({\inf\limits_{j\in\{1,2\}}}\lambda_{j}(U))\;\text{and}
\;c_{2}={\sup\limits_{U=U_L,U_R}}({\sup\limits_{j\in\{1,2\}}}\lambda_{j}(U)), 
 \end{equation*}
where $\lambda_1(U)=u-\sqrt{gh}$ and $\lambda_2(U)=u+\sqrt{gh}$ are the eigenvalues of SW1D. In practice, we use a CFL condition
 $n_{CFL}=0.5$ at second order and $n_{CFL}=1$ at first order, with
\begin{equation}
 \Delta t\leq n_{CFL} \dfrac{\Delta x}{\max\limits_{i\in\{1,...,J\}}\left(|u_i|+\sqrt{gh_i}\right)},\label{cfl}
\end{equation}
where $J$ is the number of space cells. At second order, variables $(h_i,u_i)$ in \eqref{cfl} are replaced by the reconstructed values
 $(h_{i+1/2-},u_{i+1/2-})$ and $(h_{i+1/2+},u_{i+1/2+})$ (detailed in next section).

\subsection{MUSCL-reconstruction}\label{sec:muscl}
We define the MUSCL reconstruction of a scalar function \(s\in\mathbb{R}\) (Monotonic Upwind Scheme for Conservation Law, see \cite{vanLeer79}) by
\begin{equation}
 s_{i-1/2+}=s_i-\dfrac{\Delta x}{2}Ds_i\;\text{and}\;s_{i+1/2-}=s_i+\dfrac{\Delta x}{2}Ds_i
\end{equation}
with the operator
\begin{equation}
 Ds_i=\text{minmod}\left(\dfrac{s_i-s_{i-1}}{\Delta x},\dfrac{s_{i+1}-s_i}{\Delta x}\right)
\end{equation}
and the minmod limiter
\begin{equation}
 \text{minmod}(x,y)=\left\{\begin{array}{ll}
                      \min(x,y) & \text{if}\; x,y\geq 0\\
		      \max(x,y) & \text{if}\; x,y\leq 0\\
		      0 & \text{else}
                     \end{array}\right.
\end{equation}
As mentioned previously, the MUSCL reconstruction is performed on $u$, $h$ and $h+z$, then we deduce the reconstruction of $z$. In order to keep
 the discharge conservation, the reconstruction of the velocity has to be modified as what follows \cite{Bouchut04}
\begin{equation*}
 u_{i-1/2+}=u_i - \dfrac{h_{i+1/2-}}{h_i}\dfrac{\Delta x}{2} D u_i
\quad\text{and}\quad 
 u_{i+1/2-}=u_i + \dfrac{h_{i-1/2+}}{h_i}\dfrac{\Delta x}{2} D u_i.
\end{equation*}
Thus, we recover the following conservation properties
\begin{equation*}
 \dfrac{h_{i-1/2+}+h_{i+1/2-}}{2}=h_i
\quad\text{and}\quad
 \dfrac{h_{i-1/2+}u_{i-1/2+}+h_{i+1/2-}u_{i+1/2-}}{2}=h_i u_i.
\end{equation*}

\subsection{Green-Ampt infiltration model}\label{GreenAmpt}

Infiltration is computed at each cell by using a modified Green-Ampt model \cite{Esteves00}. The soil water movement is assumed to be in the form
 of an advancing front (located at $Z_f^n$ [m]) that separates a zone still at the initial soil moisture $\theta_i$ from a saturated zone
 with soil moisture $\theta_s$ (as illustrated on figure \ref{SW2D-fig}-b. At time $t=t_n$, infiltration capacity, $I_C^n$ [m/s], is calculated
 thanks to
\begin{equation*}
 I_C^n=K_s\left(A+\dfrac{h_f-h_{over}^n}{Z_f^n}\right)\quad\text{where}\quad Z_f^n=\dfrac{V_{inf}^n}{\theta_s - \theta_i},
\end{equation*}
where $h_f$ is the wetting front capillary pressure head, $K_s$ the hydraulic conductivity at saturation, $h_{over}^n$ the overland flow
 water height (obtained from the SW2D) and $V_{inf}^n$ the infiltrated water volume. Thus we have the infiltration rate (necessary
 to couple SW2D with Green-Ampt model)
\begin{equation*}
 I^n=\dfrac{\min(h_{over}^n,\Delta t. I_C^n)}{\Delta t},
\end{equation*}
and the infiltrated volume
\begin{equation*}
 V_{inf}^{n+1}=V_{inf}^n+\Delta t.I^n,
\end{equation*}
where $\Delta t$ is the time step fixed by the CFL stability condition (see section \ref{sec:flux}).

\section{FullSWOF\_2D}
\label{sect:fullswof}
\subsection{The software}

The name FullSWOF\_2D stands for ``Full Shallow Water equations for Overland Flow in two dimensions of space''.
 It is a C++ code (free open source software under the GPL-compatible license CeCILL-V2.
 Sources can be downloaded from \url{http://www.univ-orleans.fr/mapmo/soft/FullSWOF/}) developed in the context of the project
 ANR METHODE (see \cite{Delestre10b} and \url{http://www.univ-orleans.fr/mapmo/methode/}). FullSWOF\_2D is based on a finite volume
 method (described in section \ref{sect:method}) on a structured mesh in two space dimensions. Structured grids have been chosen because
 on the one hand digital topographic maps are often provided on such meshes, and, on the other hand, it allows to develop numerical schemes in one
 space dimension (implemented in FullSWOF\_1D), extension to 2D being straightforward. As argued previously, finite volume scheme ensures by
 construction the conservation of the water mass, and is coupled with the hydrostatic reconstruction \cite{Audusse04c,Bouchut04} to deal with the
 topography source term. It preserves water height positivity and it is well-balanced ({\it i.e.} it preserves at least hydrostatic equilibrium:
 lakes and puddles). Several numerical fluxes (Rusanov, HLL, kinetic and VFRoe-ncv flux) and second order reconstructions (MUSCL, ENO and
 modified ENO reconstruction) are implemented. Currently, we recommend, based on \cite{Delestre10b}, to use the second order scheme with MUSCL
 reconstruction \cite{vanLeer79} and HLL flux \cite{Harten83} detailed in section \ref{sect:method}. FullSWOF\_2D is designed
 in order to ease implementation of new numerical methods. FullSWOF\_2D has already been validated on analytical solutions integrated
 in SWASHES (Shallow Water Analytic Solutions for Hydraulic and Environmental Studies: a free library of analytical solutions \cite{Delestre2012})
 and on real events simulation at small scales \cite{Delestre10b,DBLP:journals/corr/abs-1204-3210} (on parcels). In next section,
 FullSWOF\_2D is validated at a bigger scale on the Malpasset test case.

\subsection{Validation: Malpasset case}

 \begin{figure}[htbp]
\begin{center}
\includegraphics[width=12cm,height=5.5cm]{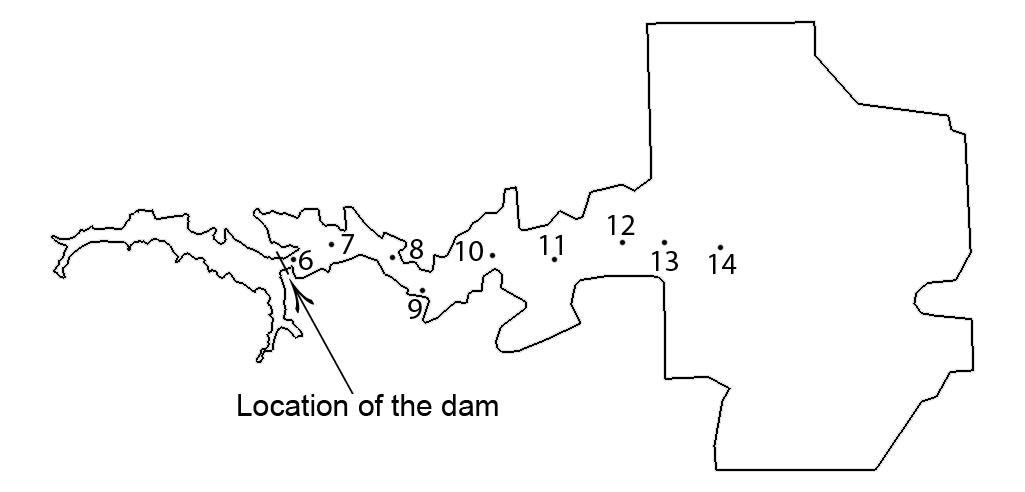}
\caption{\label{ffglut} Position of the different gauges and location of the dam.}
\end{center}
\end{figure}

Malpasset dam break is a real life application. It occured (because of heavy rains) on December 2 1959 in the Var Department in southeastern
 France, causing 433 casualties. This dramatic event is very often mentioned in studies dealing with dam breaks and associated
 risks (see among others \cite{Goubet79,Lebreton85,Benoist89,Gioda02,Chanson05a,Chanson06,Mulder09}). Because of its varying topography and complex
 geometry, it is widely used as a benchmark test for numerical methods and hydraulics software
 \cite{Hervouet99,Hervouet00,Malleron11,Singh11,Audusse05,Duran12,Ata}. Moreover, this problem allows to test the ability of the scheme to treat
 the still water (at the level of the sea downstream, before the wave reaches it) and the wet-dry interfaces.
 
The dimensions of the domain are
$dim_x=17273.9\;\text{m}$ along x-axis with 1000 meshes and $dim_y=8381.3\;\text{m}$ along y-axis with 486 meshes, the total time of simulation
 is $2500\;\text{s}$ and we consider the Manning law with $n=0.033$ as advised in the literature \cite{Hervouet00}.

The dam is considered as a straight line (see figure \ref{ffglut}), the water level inside the reservoir is set to 100 m above sea level and
 the computational domain downstream the dam is considered as dry bed. Indeed, the initial discharge in the river before dam failure can be
 neglected because of the huge amount of flow caused by the dam failure. 

\begin{figure}[!htb]
\begin{center}
\includegraphics[width=12cm,height=5.5cm]{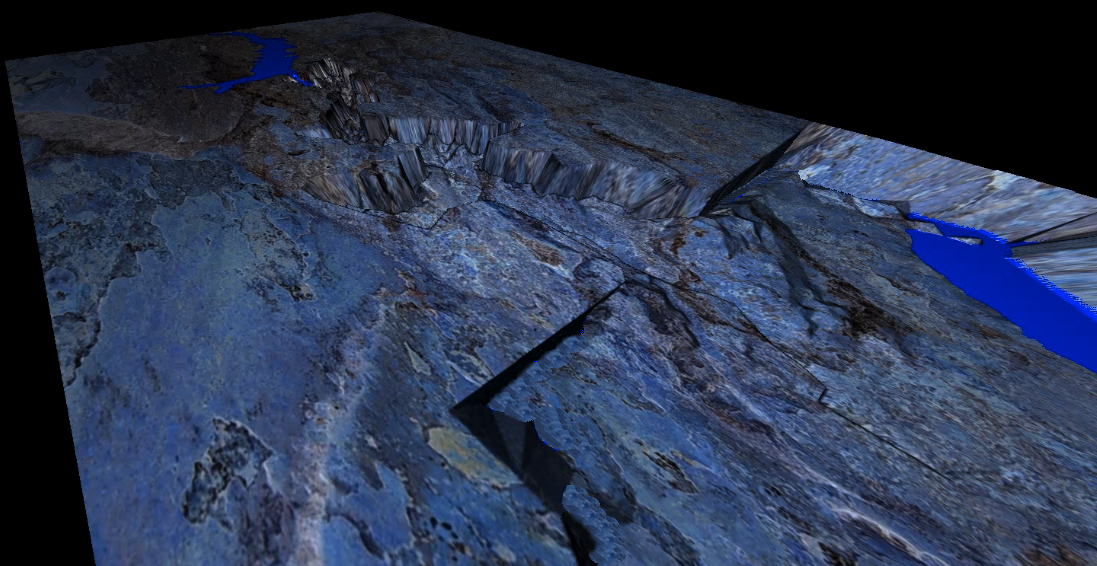}

\vspace{.5cm}

\includegraphics[width=12cm,height=5.5cm]{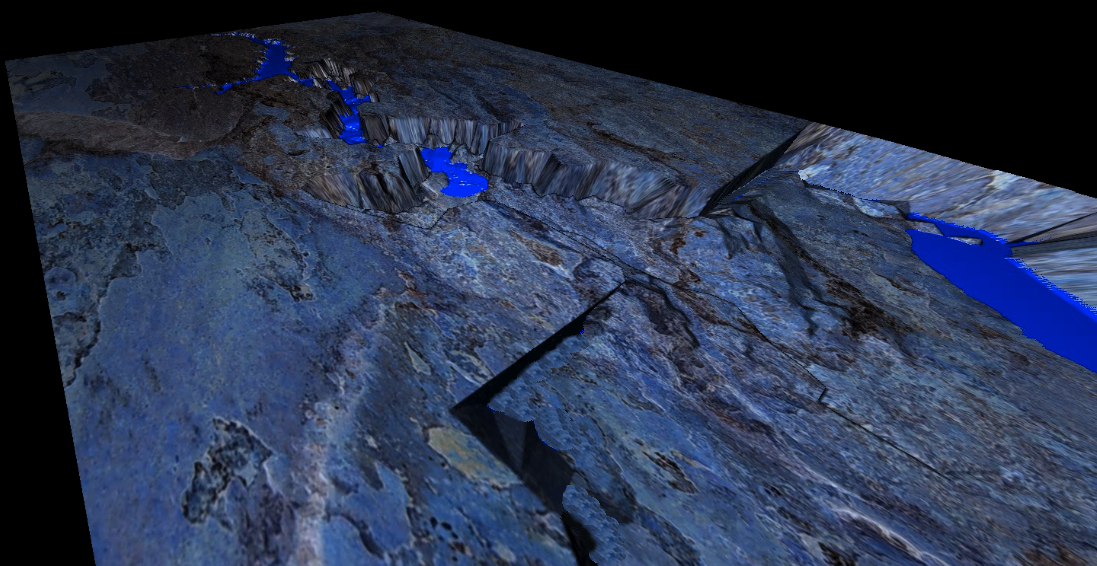}

\caption{\label{malpasset_0}Propagation of the Malpasset dam break for different time 0 s and 500 s.}
\end{center}
\end{figure}

\begin{figure}[!htb]
\begin{center}

\includegraphics[width=12cm,height=5.5cm]{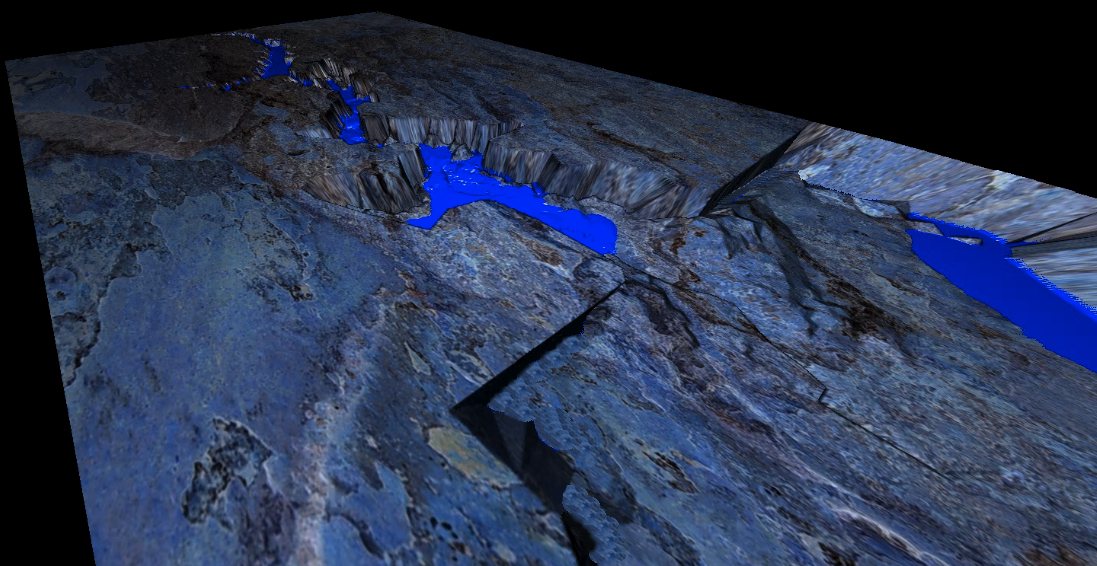}
\vspace{.5cm}

\includegraphics[width=12cm,height=5.5cm]{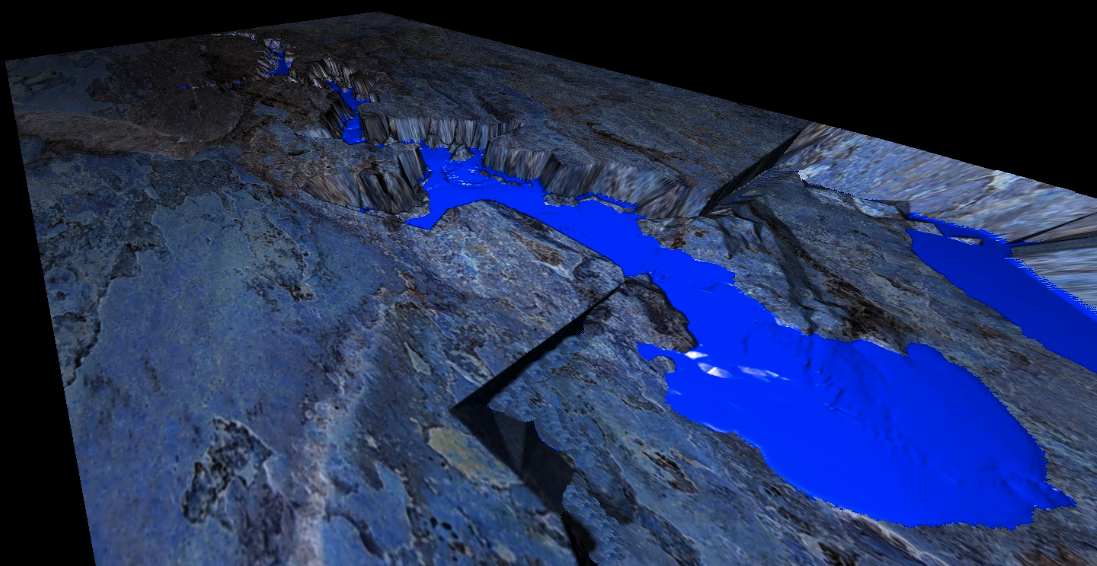}

\vspace{.5cm}

\includegraphics[width=12cm,height=5.5cm]{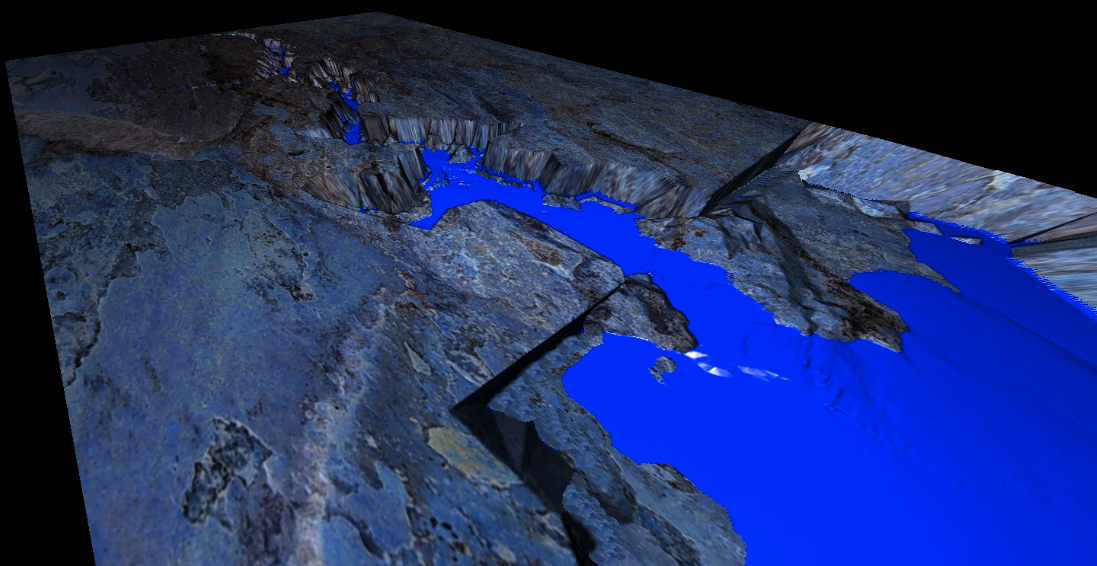}

\caption{\label{malpasset_50}Propagation of the Malpasset dam break for different time 1000s, 2000 s and 2500 s.}
\end{center}
\end{figure}

We run the simulation using FullSWOF\_Paral over 16 processors and we get the
 results in Figures \ref{malpasset_0} and \ref{malpasset_50}. 
In 1964, a physical model wih a scale of 1:400 was built by Laboratoire National
 d'Hydraulique to study the dam-break flow. The maximum water level and flood wave arrival time were recorded at 9 points in the physical model
 (named from 6 to 14). This physical model has been used to calibrate/validate Telemac\_2D software \cite{Hervouet99,Hervouet00,Hervouet07}
 and these results are used here to validate FullSWOF\_2D at ``big'' scale. We can see on figure 5 that the results are very closed to those
 obtained with Telemac. On figures \ref{malpasset_0} and \ref{malpasset_50} is represented the propagation of the wave due to the dam break 
 (a video can be visualized on \url{http://www.youtube.com/user/FullSWOF}) and on figure 6, the propagation of the water height at the gauges during time (this will be used in
 future work). These results show the efficiency of FullSWOF but at that scale, details such as houses and roads are not represented. So in case
 of prevention, the use of these results might be limited. It might be useful to have more details. So the mesh would be finer, in that case the
 code needs to be parallelized.

\begin{figure}[!h]
\begin{center}
\resizebox{8cm}{!}{\includegraphics{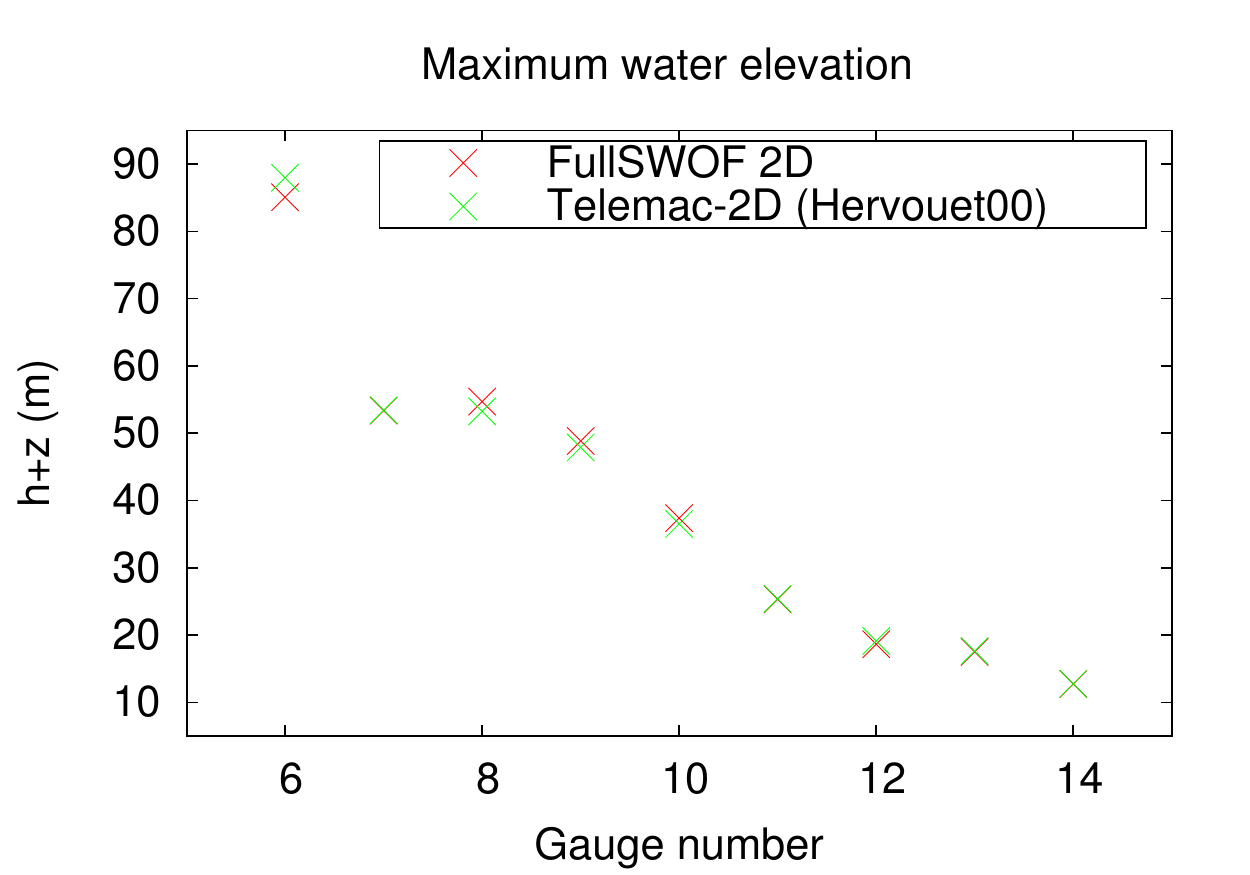}}
\caption{Comparison of the maximum height of the water level in all gauges between FullSWOF\_2D and Telemac\_2D.}
\end{center}
\end{figure}

\begin{figure}[!h]
\begin{center}
\resizebox{8cm}{!}{\includegraphics{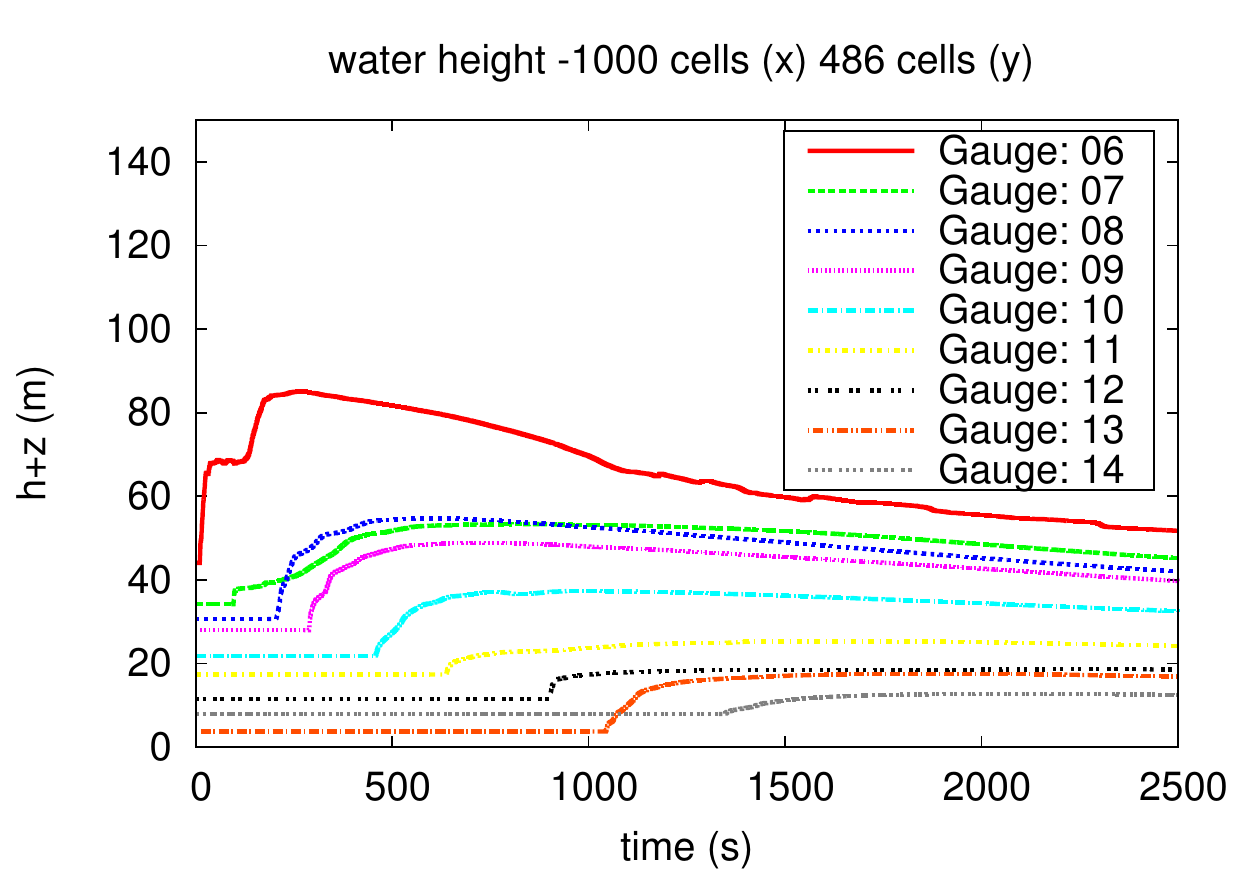}}
\caption{Propagation of the water height for different gauges during time.}
\end{center}
\end{figure}

\section{SkelGIS}
\label{sect:skelgis}
\subsection{Purpose}
Parallelism~\cite{Kumar:2002:IPC:600009} is an intricate domain of computing science. Indeed, in addition to good sequential and parallel
 programming skills, it requires also a strong knowledge on processors, memory and network use to be able to write optimal programs on modern
 parallel computers. Writing a parallel program is long and complex, and implementing optimizations to get very good performances is even more
 difficult.
As a consequence, giving the opportunity to non-specialists and even to non-computer scientists to use parallelism is a major research area for a
 long time. Indeed, in lots of domains high performance computing has become a crucial point.

A wide range of solutions exists to give access to parallelism to non-specialists. Most of them require the collaboration of experts in parallelism
 with the scientists of the targeted domain since classical parallel libraries are too technical. A good alternative is to propose libraries that
 clearly separate the parallelism side from the user side.
Those tools provide some parallel patterns that hide the technical details of parallelism and are optimized for a class of problems. Such libraries
 give a restricted access to parallelism but greatly simplify the program development. The difficulty is to find the good balance between the level
 of abstraction, to hide the parallelism, and the performance of the program.

SkelGIS library is presented in this part and takes place in the latter category of libraries. SkelGIS is an algorithmic skeleton library and has
 been used to parallelize FullSWOF\_2D.

\subsection{Algorithmic skeletons}

Algorithmic parallel skeletons were introduced in 1988 by Muray Cole \cite{Cole1}. The domain of algorithmic skeletons aims at providing generalist
 patterns of parallelization to the user and hide any parallel implementation of the resulting application in the library. As a result only
 sequential interfaces of the library (named skeletons) are used by scientists and a parallel application is executed without any knowledge on
 parallel libraries such as MPI, OpenMP, CUDA etc. We can enumerate lots of algorithmic skeleton libraries, all of them based on Muray Cole's
 work \cite{Cole1,Benoit:2005:TFC:2152883.2152988,Cole:2004:BSO:1007980.1007985}. Some of them are implemented in Java, some of them in C++ as for
 example QUAFF~\cite{DBLP:journals/pc/FalcouSCL06}, eSkel~\cite{DBLP:conf/europar/BenoitCGH05}, Muesli~\cite{DBLP:conf/europar/BotorogK96},
 SkeTo~\cite{DBLP:conf/infoscale/MatsuzakiIEH06} and OSL~\cite{DBLP:conf/pact/JavedL11}.

To explain the concept of algorithmic skeleton, the equation \eqref{map} represents the well known \emph{map} skeleton where $F$ is the set of
 functions expressed by equation \eqref{func} and DStruct the type of elements in the distributed data structure given to the map skeleton. A
 \emph{map} skeleton takes an input data structure, and a user function in inputs, apply the user function to each element of the data structure
 and return a new resulting data structure. Parallelism is hidden in the data structure that is distributed transparently, and in the skeleton
 call. To get a parallel code with a map, the user simply has to construct a data structure $d$ and to define a sequential function $f$
 \eqref{func} that gives the calculation to accomplish on one element of the data structure. Then he or she has to call the skeleton this
 way : map($d$, $f$).

\begin{equation}\label{map}
\emph{map}: F \times DStruct \longrightarrow DStruct
\end{equation}
\begin{equation}\label{func}
F = \left\{ f:E \longrightarrow E \right\}
\end{equation}

Though, this type of skeletons are not enough expressive. For example, a classical operation on a two dimensional dataset consists in computing
 a value in a matrix from the four or eight neighbor values around this point. With skeletons of type \emph{map}, it is not possible to directly
 access this neighborhood since only the current element is available. To sidestep this limitation existing libraries propose communication
 skeletons, and especially the \emph{shift} skeleton that permits to shift the matrix in order to get access to another value in the matrix.
 This skeleton has two major issues. First, the shifted matrix implies a copy of the matrix which is undesirable when dealing with a huge
 amount of data. Second, the shift skeleton makes difficult the programming of classical algorithms since it must be
 used in place of a classical matrix element accessors.

\begin{figure}[!h]
\begin{center}
\resizebox{8cm}{!}{\includegraphics{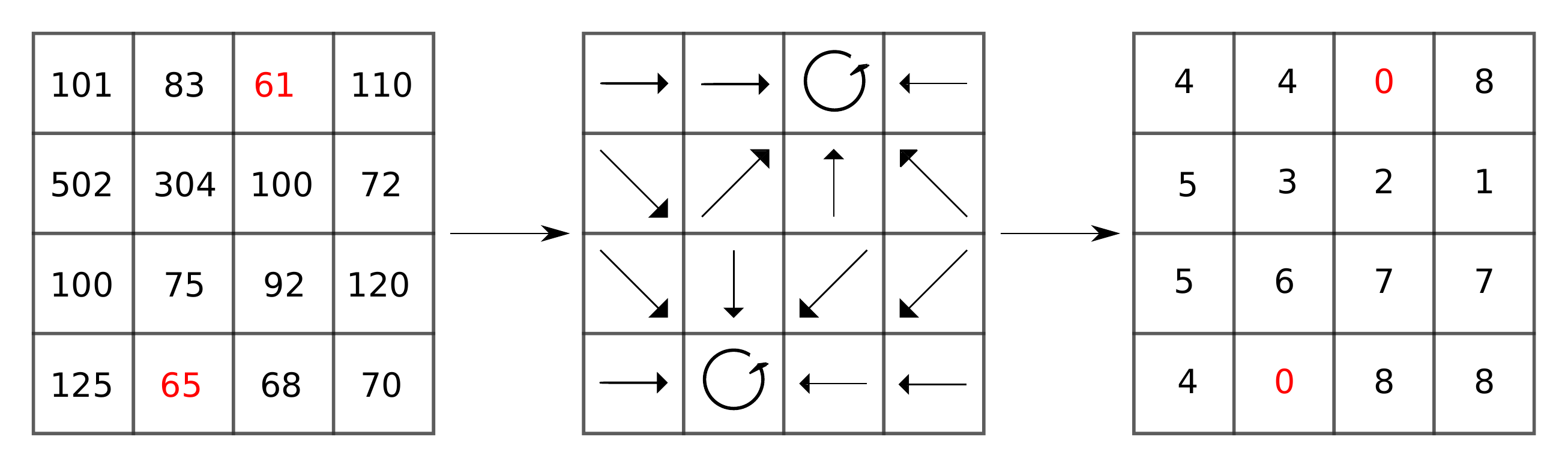}}
\caption{Flow direction algorithm}
\label{directions}
\end{center}
\end{figure}

To illustrate this, an example of an eight neighbors algorithm is used. Geo-scientists are working on the automatic delination of watershed to
 determine, for example where the water is going to flow, to elaborate irrigation systems, or to determine the impact of a water flooding on a
 land. These calculations are performed on topography matrices. One important and simple step of this calculation is the flow direction
 computation (figure \ref{directions}). It consists in determining for each point of the DEM where the water would flow with the assumption that
 it should flow in the direction of the nearest lower point. Therefore, the flow direction algorithm is simple, it consists in getting the minimum
 value for the eight neighbors around the current point and saving the corresponding direction. 

Eight neighbors algorithms are a good example of algorithm used to solve scientific problems. Indeed, it is a frequent algorithm in mathematics
 (finite volumes discretization), physics, image processing, geo-sciences and many other domains. However, it is not easy to make the flow
 direction calculation parallel with standard skeletons (map, shift etc.). Figure \ref{sketo} shows the 24 skeleton calls needed. Height
 \emph{shift} skeletons are needed to get height neighbors values, then nested calls of \emph{zipwith} skeletons are needed to get the minimum
 value, finally \emph{map} skeletons are needed to associate the minimum value to a direction.

\begin{figure}[h!]
\small
  \begin{center}
    \begin{lstlisting}[mathescape,basicstyle=\footnotesize,frame=single,language=C++]
zipwith(min,zipwith(min,zipwith(min,zipwith(min,map(f,m),
map(f,shift<1,1>(shift_id,m))),zipwith(min,
map(f,shift<1,0>(shift_id,m)),map(f,shift<1,-1>(shift_id,m)))),
zipwith(min,zipwith(min,map(f,shift<0,-1>(shift_id,m)),
map(f,shift<-1,-1>(shift_id,m))),zipwith(min,
map(f,shift<-1,0>(shift_id,m)),map(f,shift<-1,1>(shift_id,m))))),
map(f,shift<0,1>(shift_id,m)));
\end{lstlisting}
\caption{24 nested calls needed for the flow direction calculation with standard skeletons}
\label{sketo}
\end{center}
\end{figure}

Then, this type of skeleton libraries have few problems :
\begin{itemize}
\item The difficulty of getting a parallel code is moved to the difficulty of functional programming with skeletons.
\item The code is difficult to write and even more difficult to read.
\item Performance damages occur because of all the skeleton calls and because of the duplication of matrices with the shift skeleton.
\end{itemize}

\subsection{SkelGIS}

SkelGIS is a new kind of algorithmic skeleton library. The objective is to propose a hierarchical skeleton library for very large matrices and
 to give access to parallelism through sequential adapted interfaces in order to keep a sequential programming style. SkelGIS gets out of the
 functional paradigms in order to stick to the programming habits of non-specialists. It relies on two concepts and a distributed data structure.

\subsubsection{Concepts}

\emph{The first concept of SkelGIS is to propose basic skeletons where the user function has a direct access to a matrix. The user has to define
 functions of type $f:M \longmapsto M$ ($M$ is a SkelGIS matrix).}
\medskip

This is a major difference with the classical skeletons (as map) where the user functions defines what to do on a single matrix element
 (functions of type $f:E \longmapsto E$ where $E$ is the type of elements of a matrix $m$). The user describes, in the sequential function, what
 to do on the matrix. As a result, the expressiveness of SkelGIS is very good and no communication skeletons are needed. SkelGIS offers a way to
 stay close to a sequential programming style. Another consequence of the expressiveness of basic skeletons is that a single call of skeleton at
 a time is needed, and performances of SkelGIS are better than existing skeletons. For example, the figure \ref{perf} gives the performances
 of the flow direction calculation, above-cited, with SkelGIS and SkeTo.

\begin{figure}[!h]
\begin{center}
\begin{tabular}{|c|c|c|c|}
  \hline
  \textbf{Size}  & \textbf{SkeTo} & \textbf{SkelGIS}\\
  \hline
  \textbf{$339 \times 225$} & 30ms & 21ms\\
  \hline
  \textbf{$14786 \times 10086$} & 101s & 31s\\
  \hline
\end{tabular}
\caption{Comparison between SkeTo library and SkelGIS, on flow direction calculation}
\label{perf}
\end{center}
\end{figure}

\begin{figure}[!h]
\begin{center}
\resizebox{9cm}{!}{\includegraphics{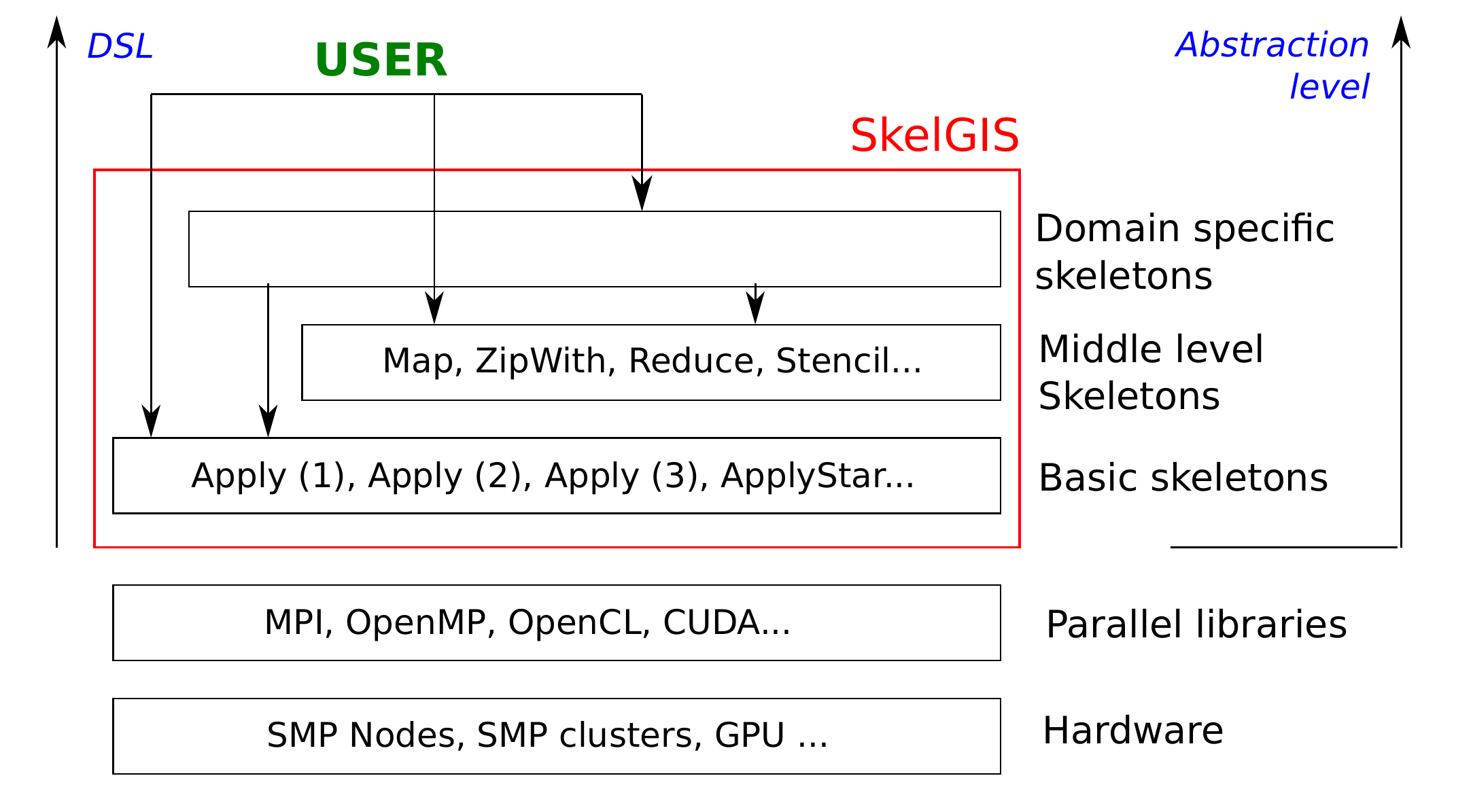}}
\caption{Hierarchy of SkelGIS}
\label{hie}
\end{center}
\end{figure}

\medskip
\emph{The second concept of SkelGIS is to propose a hierarchy of skeletons (Figure~\ref{hie}) where every higher abstraction level skeleton uses
 basic skeletons.}

As a result, SkelGIS skeletons inherit optimizations and hardware support of basic skeletons. New optimizations or new hardware support only has
 to be added in basic skeletons to be available in the whole library. This hierarchy also provide a clear abstraction choice for the user.
Ensuring the scalability and durability of a library is an important feature to make it live and used. Skeleton libraries are based on parallel
 libraries depending on hardware, and hardware can completely change quite quickly. Existing libraries propose a set of skeletons that are
 independent from each others. As a result, taking into account a new hardware (or a new parallel library) requires to re-implement all the
 skeletons of the library.

\subsubsection{Distributed data structure}

As every skeleton library, SkelGIS skeletons are applied to a distributed data structure. The current version of SkelGIS only proposes a
 distributed two-dimensional matrix to apply skeletons on. This distributed matrix named \emph{DMatrix} is the basis of every parallel
 calculation made by SkelGIS. The constructor of the DMatrix is responsible for dividing the initial matrix in blocks and distributing the load
 of each block from the disk on different processors.

This DMatrix is used by skeletons to apply the user functions in parallel on the whole matrix. As the user can manipulate distributed matrices
 without being aware of the parallel technical details, some basic tools have been developed to manipulate them. Then, manipulation of the
 \emph{DMatrix} looks like data structure of the \emph{Standard Template Library} (STL in C++) manipulation. The three main tools to manipulate
 a \emph{DMatrix} are:
\medskip
\begin{itemize}
\item A set of iterators to navigate in the distributed matrix.
\item A set of neighborhood functions to get the neighbors of the current element.
\item Get/Set to get a value of the distributed matrix or to write a value in the distributed matrix.
\end{itemize}
\medskip

The user uses the DMatrix as if it was a common sequential matrix of the STL, and the skeletons are in charge of hiding block division and
 communication aspects. Figures \ref{main} and \ref{code} represent the user code needed to get the parallel version of the flow direction
 calculation above-cited. The main function is in charge of the initialization of the library SkelGIS, the construction of matrices and the call
 of skeletons. Then, user functions have to be programmed, in a sequential programming style, using the \emph{DMatrix} tools.

\begin{figure}[!h]
\begin{center}
\begin{lstlisting}[mathescape,frame=single,language=C++]
#include ``skelgis.h''
int main(int argc, char** argv)
{
    //initialization of the library
    INITSKELGIS;
    //creation of matrix m
    DMatrix<float> * m = new DMatrix<float>(file,1);
    //call of the skeleton apply with the function
    // "directions" on the matrix m
    DMatrix<int> * m2 = apply<float,int>(directions,m);
    //end of the library use
    ENDSKELGIS;
}
\end{lstlisting}
\caption{Main function of the user program}
\label{main}
\end{center}
\end{figure}

\begin{figure}[!h]
\begin{center}
\begin{lstlisting}[mathescape,frame=single,language=C++]
//begin declaration of the function directions
BEGINApply(directions, input, float, output, int)
{
  //initilization of iterators on input
  //and output matrices
  DMatrix<float>::iterator it = input->begin();
  DMatrix<int>::iterator itOut = output->begin();
  DMatrix<float>::iterator itEnd = input->end();
  //8 neighbors table
  DMatrix<float>::neighbors nghb;
  //for each element of the matrix
  for(; it <= itEnd; ++it , ++itOut)
    { 
      int index=0; 
      //get the current value of the iterator
      float min = input->getValue(it);
      //get the 8 neighbor values
      input->get8Neighbors(it,nghb);
      //get the minimum
      for(int i=0;i<8;i++){
          if(nghb[i]>0 && nghb[i]<min) {
              index=i+1;
              min=nghb[i];}
      }
      //set the minimum direction in the output matrix
      output->setValue(index, itOut);
    }
}
//end declaration of the function
END(directions);
\end{lstlisting}
\caption{User function given to the skeleton}
\label{code}
\end{center}
\end{figure}

Behind SkelGIS, the exact same work as in MPI is done, but automatically. Actually, the library proceeds to a domain decomposition of a
 structured matrix and manages MPI exchanges of ghost cells at the border of the domain decomposition. Then SkelGIS proposes an easy way to
 code structured approaches of simulations. Of course, the structured approach is limited regarding the complexity of some simulations, but
 SkelGIS is under development and will propose in the future unstructured distributed data structures and also graph and tree distributed data
 structures. This way, SkelGIS will be able to manage different types of simulations.

\section{Parallelization of FullSWOF\_2D}
\label{sect:paral}
\subsection{MPI}
In this section, we will explain the domain decomposition method which we have used to parallelize FullSWOF\_2D (for more details see
 \cite{Brugeas96}). This method has been applied by using the library of functions MPI (Message Passing Interface). This method has been chosen
 in order to be the least intrusive as possible in the FullSWOF structure to ease future development. In fact, the domain decomposition method
 consists in splitting the  data into many independent domains. On each domain, the code is executed by one anonymous process and when it is
 necessary the process communicates via calls to MPI communication primitives. MPI is a standardized and portable message-passing system used to
 compute on both shared memory and distributed memory machines. The parallelism with MPI decomposes into four main steps:
\begin{enumerate}
\item decomposition of the domain,
\item knowing his four neighbors,
\item exchanging the points on the interfaces
\item and calculating the scheme on each process.
\end{enumerate}
From these four steps, we have parallelized FullSWOF\_2D.

At the beginning of the run, we create a group of $n$ processes. Each process is assigned a rank between $0$ and $n-1$. For us, this linear
 ranking of processes is not appropriate, because the structure of the computational domain is in two dimensions with communications into the
 two directions. So, the processes are arranged in the 2D Cartesian topology thanks to MPI topology mechanism.

 MPI creates a relationship between ranks and Cartesian coordinates and we used these information to know the neighbors of each nodes in order to
 communicate the values of the interfaces. After this step, the domain can be represented by a graph where the nodes stand for the processes
 and the edges connect process that communicate with each other. So, we can compute the number of points for each node in $x$ and $y$ and to
 initialize the local variables.

After the initialization, we exchange the value of the interfaces to compute the value of the variables at the next step of the time. Indeed,
 when we compute the value of the variables according to the order of the scheme. At order one, we add a line and column of the cells that
 receive the value necessary to the computation of the scheme and if we are at order two, we add two lines and two columns.

Concerning the communication between the nodes, MPI provides message-passing between any pair of processes. So, having exchanged the values at
 the interfaces, we compute the next step of the time by applying the algorithm of sequential FullSWOF. 

\subsection{SkelGIS}

Using SkelGIS to parallelize the shallow water equations is very different from a standard parallelization. No parallel knowledge
 (domain decomposition, communications etc) and no parallel code (MPI use) are needed at all to get a parallel version of the simulation.
 However, the simulation has to be thought with skeleton use. 

Algorithm \ref{algo} represents the pseudo code of the shallow water equations simulation and, in red is added the interventions needed to
 use SkelGIS with this algorithm. First, the initialization of matrices has to be made with the \emph{DMatrix} object, then each algorithm step
 has to be executed through a skeleton call (\emph{applyList skeleton}). Finally, the sequential code has to be a little modified because accesses
 to matrices elements have to be done through \emph{DMatrix} iterators and accessors.

\begin{algorithm}[!h]
 Initialization of variables \textcolor{red}{with DMatrix}\\
 \For{$t\leftarrow 0$ \KwTo $N$}{
  \For{$i\leftarrow 0$ \KwTo $2$}{
    \textcolor{red}{applyList : }Boundary Conditions\\
    \textcolor{red}{applyList : }Second order reconstruction\\
    \textcolor{red}{applyList : }Hydrostatic reconstruction\\
    \textcolor{red}{applyList : }Numerical flux computation\\
    \textcolor{red}{applyList : }Scheme computation\\
     \textcolor{red}{applyList : }Frictions
   }
   \textcolor{red}{applyList : }Order 2 along $t$ - Method of Heun
 }
 \caption{Shallow water equations algorithm with additional calls in red for SkelGIS}
 \label{algo}
\end{algorithm}

In addition to the fact that no parallelization knowledge is needed using SkelGIS, the non-computing scientist is fully independent in parallel
 developments. The parallel application getting from the use of SkelGIS is organized and easily upgradeable. Furthermore, an OpenMP version and
 the use of SSE instructions (Streaming SIMD extensions) optimizations are under development in SkelGIS. As a result with no efforts the code
 of shallow water equations will work with OpenMP for shared memory systems as with Intel additional processor instructions.

\subsection{Results}

\begin{figure}[!h]
\begin{center}
\resizebox{8cm}{!}{\includegraphics{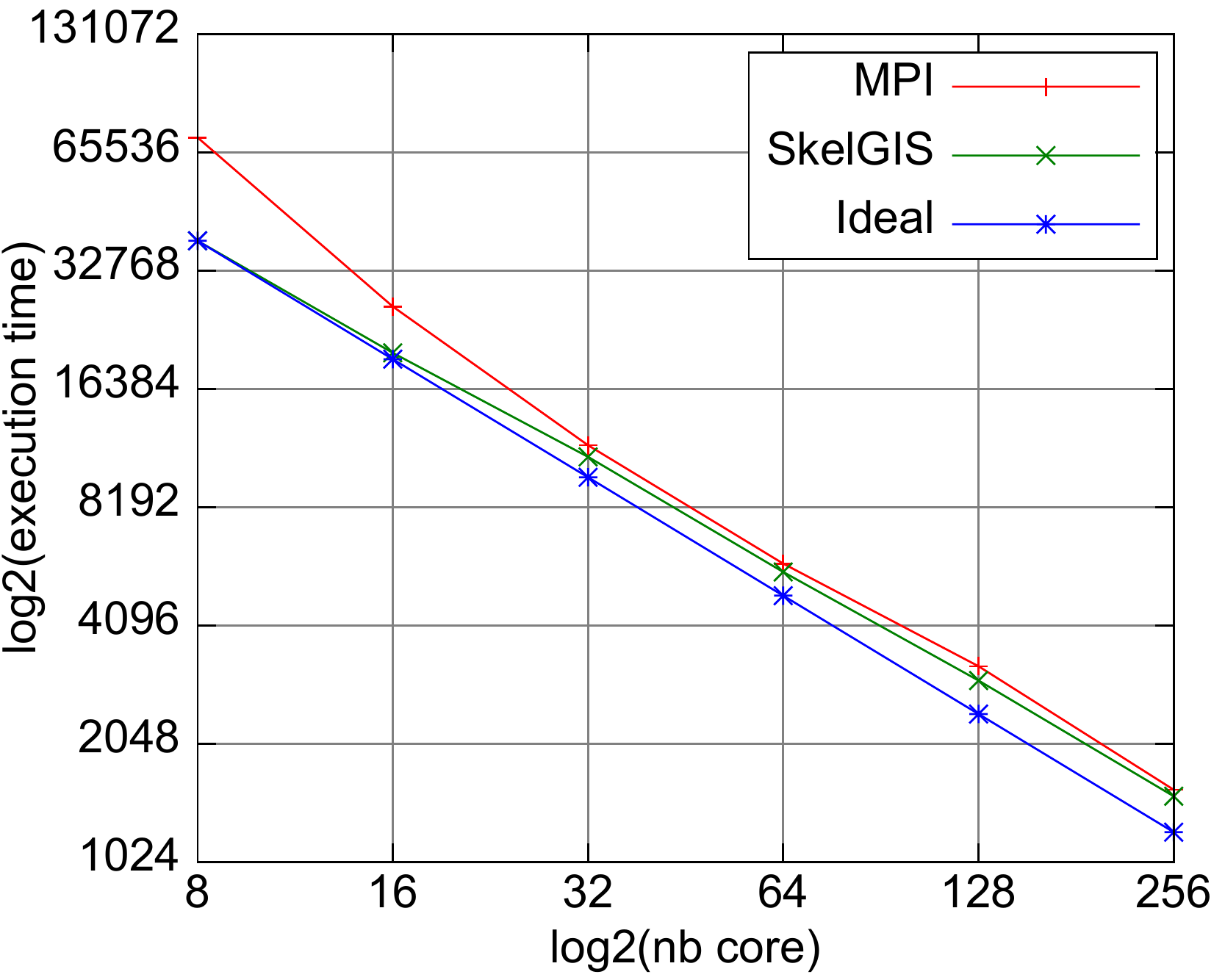}}
\caption{Comparison between MPI, SkelGIS and Ideal versions of shallow water equations for 5000 time iterations : Log 2 of execution times}
\label{5000}
\end{center}
\end{figure}

\begin{figure}[!h]
\begin{center}
\resizebox{8cm}{!}{\includegraphics{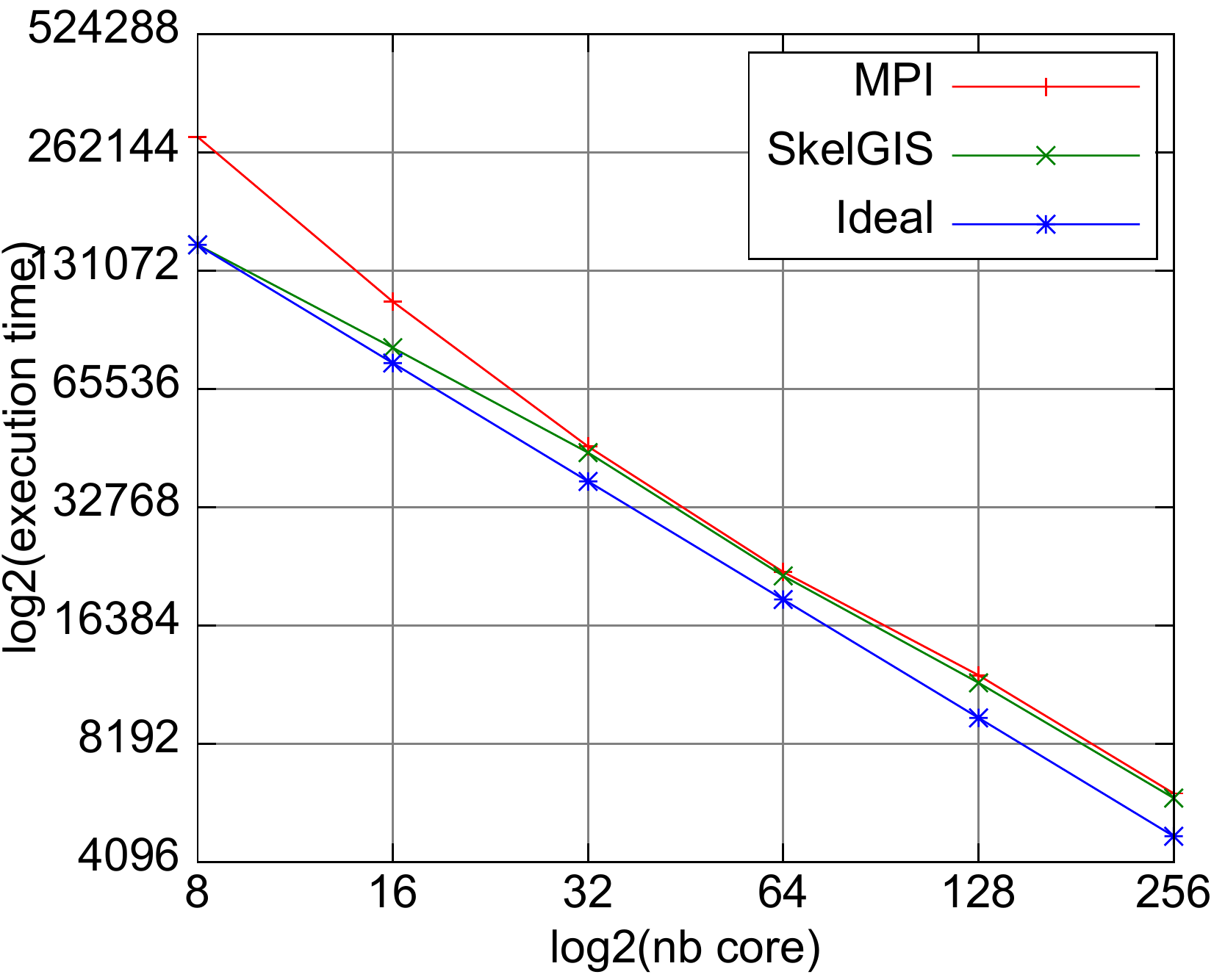}}
\caption{Comparison between MPI, SkelGIS and Ideal versions of shallow water equations for 20000 time iterations : Log 2 of execution times}
\label{20000}
\end{center}
\end{figure}

In this part are presented some comparison results between the MPI implementation and the SkelGIS implementation of the shallow water equations.
 These benchmarks were made on a 5120x5120 domain with 5000 iterations and 20.000 iterations. Figures \ref{5000} and \ref{20000} represent
 the base two logarithm of the execution times. Then, it represents the inverse of the speedup and the ideal execution time is drawn in blue.

Both implementations have very good results and are closed to the ideal execution time. We can notice, however, that SkelGIS results are linear
 while the MPI version starts with a bigger execution time with 8 processors. From 32 to 256 processors the MPI implementation meets SkelGIS
 execution times.

\subsection*{Conclusions}

As a conclusion, both implementations are efficient, but parallelization concepts of these two implementations are very different. The Skeletton
 strategy requires a more important work to rewrite the code, using SKELGIS library but it provides a final code with a greater expressivity
 and the parallelization technics are hidden within the library. On the contrary, the MPI version involves less changes in the software, but it
 requires to deal with the parallelization concept. Therefore, our first conclusion is that the performances of the two approaches are very
 close. For an existing code, depending on the skills of the developers, each strategy can be preferred. If the choice is done at the very
 beginning of the project, skeletons offers a more expressive code and does not require to know about parallelization. Both FullSWOF\_2D and
 SKELGIS are still under development, interested readers might contact us for further informations.


\end{document}